\documentclass{amsproc}

\usepackage{latexsym}
\usepackage{amsmath}
\usepackage{amsfonts}
\usepackage{amssymb}
\usepackage{amscd}

\numberwithin{equation}{section}

\newtheorem{proposition}{Proposition}[section]
\newtheorem{lemma}[proposition]{Lemma}
\newtheorem{theorem}[proposition]{Theorem}
\newtheorem{corollary}[proposition]{Corollary}
\newtheorem{definition}[proposition]{Definition}
\newtheorem{example}[proposition]{Example}

\def\itheorem#1#2{\newtheorem{#1}[proposition]{#2}}

\itheorem{Remark}{Remark}

\def\Hom{\mathop{\rm Hom}\nolimits}
\def\Ext{\mathop{\rm Ext}\nolimits}

\def\rd{\mathop{\rm rd}\nolimits}
\def\im{\mathop{\rm Im}\nolimits}
\def\dom{\mathop{\rm Dom}\nolimits}

\def\Q{{\mathbb{Q}}}
\def\Z{{\mathbb{Z}}}

\begin{document}

\title{Filtration-equivalent $\aleph_1$-separable abelian groups of cardinality ~$\aleph_1$}

\author{Saharon Shelah}

\address{Department of Mathematics,
The Hebrew University of Jerusalem, Israel, and Rutgers
University, New Brunswick, NJ U.S.A.}

\email{Shelah@math.huji.ac.il}
\thanks{
Number 855 in Shelah's list of publications. The first author was
supported by project No. I-706-54.6/2001 of the {\em
German-Israeli Foundation for Scientific Research \&
Development}.\\ The second author was supported by a grant from
the German Research Foundation DFG}

\author{Lutz Str\"ungmann}

\address{Department of Mathematics,
University of Duisburg-Essen, 45117 Essen, Germany}

\email{lutz.struengmann@uni-essen.de}

\begin{abstract} We show that it is consistent with ordinary set
theory $ZFC$ and the generalized continuum hypothesis that there
exist two $\aleph_1$-separable abelian groups of cardinality $\aleph_1$ which
are filtration-equivalent and one is a Whitehead group but the
other is not. This solves one of the open problems from
\cite{EkMe}.
\end{abstract}

\maketitle \setcounter{section}{1}

\section*{Introduction}
An $\aleph_1$-separable abelian group is an abelian group $G$ such
that every countable subgroup is contained in a free direct
summand of $G$. This property is apparently stronger than the
property of being strongly $\aleph_1$-free; however, the two
properties coincide for groups of cardinality at most $\aleph_1$
in models of Martin's Axiom ($MA$) and the negation of the
continuum hypothesis ($\neg CH$). Over the years the variety and
abundance of $\aleph_1$-separable groups obtained by various
constructions has demonstrated the failure of certain attempts to
classify $\aleph_1$-separable groups of cardinality $\aleph_1$. In
brief, one can say that positive results towards classification
can be given assuming $MA + \neg CH$ and negative results are
obtained assuming $CH$ or even G\"odel's constructible universe
$V=L$. A good survey is for instance \cite[Chapter XIII]{EkMe}.\\
There are four principal methods of constructing $\aleph_1$-free
groups: as the union of an ascending chain of countable free
groups; in terms of generators and relations; as a subgroup of a
divisible group; and as a pure subgroup of $\Z^{\omega_1}$. In the
study of $\aleph_1$-separable groups it turned out to be helpful
to consider the concept of {\it filtration-equivalence}, a
relation between two $\aleph_1$-separable groups. Recall that two
groups $A$ and $B$ of cardinality $\aleph_1$ are called {\it
filtration-equivalent} if they have filtrations $\{ A_{\nu} : \nu
\in \omega_1 \}$ and $\{ B_{\nu} : \nu \in \omega_1 \}$
respectively such that for all $\nu \in \omega_1$, there is an
isomorphism $\Phi_{\nu}: A_{\nu} \rightarrow B_{\nu}$ satisfying
$\Phi_{\nu}[A_{\mu}]=B_{\mu}$ for all $\mu \leq \nu$. Such an
isomorphism is called {\it level-preserving}. Note that it is not
required that $\Phi_{\tau}$ extends $\Phi_{\nu}$ when $\tau \geq
\nu$ and that filtration-equivalent groups $A$ and $B$ are also
quotient-equivalent, i.e. for all $\nu \in \omega_1$ we have
$A_{\nu +1}/ A_{\nu} \cong B_{\nu + 1}/B_{\nu}$.\\
Under the hypothesis of Martin's Axiom the notion of filtration-equivalence represents the end of the search; more precisely,
assuming $MA + \neg CH$, filtration-equivalent
$\aleph_1$-separable groups are isomorphic. Assuming even the
proper forcing axiom ($PFA$) every $\aleph_1$-separable group (of
cardinality $\aleph_1$) is of a special {\it standard form}.
However, in $V=L$ there exist non-isomorphic $\aleph_1$-separable
groups of cardinality
$\aleph_1$ which are filtration-equivalent.\\
In \cite[Open problems on the structure of $\Ext$ Nr.6]{EkMe}
Eklof and Mekler asked whether or not it is consistent with $ZFC$
that there exist two filtration-equivalent $\aleph_1$-separable
groups of cardinality $\aleph_1$ such that one is a Whitehead
group and the other is not. Recall that a Whitehead group is an ablian group $G$ satisfying $\Ext(G,\Z)=0$. The class of Whitehead groups is closed under direct sums and subgroups and contains the class of free abelian groups. However, the question whether all Whitehead groups are free is undecidable in ZFC as was shown by the first author in \cite{Sh1}, \cite{Sh2}. Similarly, we shall show in this paper that the answer to the question by Eklof and Mekler is affirmative even assuming $GCH$.\\

All groups are abelian and notation is in accordance with
\cite{Fu} and \cite{EkMe}. For further details on
$\aleph_1$-separable groups and set theory we refer to
\cite{EkMe}.

\section*{The Construction}
Using special ladder systems we construct $\aleph_1$-separable abelian groups
of cardinality $\aleph_1$ with a prescribed $\Gamma$-invariant.
The construction is similar to the one given
in \cite[Chapter XIII, Section 0]{EkMe}.\\

Throughout this paper let $S \subseteq \omega_1$ be a stationary
and co-stationary subset of $\omega_1$. Since $\lim(\omega_1)$ is
a closed and unbounded subset of $\omega_1$ we may assume without
loss of generality that $S$ consists of limit ordinals of
cofinality $\omega$ only. We shall further require that $\omega^2$
divides $\delta$ for every $\delta \in S$. We recall the
definition of a ladder and a ladder system, respectively.

\begin{definition}We use the following notions:
\begin{enumerate}
\item  A {\rm ladder} on $\delta \in S$ is a strictly increasing
sequence $\eta_{\delta}=\{\eta_{\delta}(n): n \in \omega \}$ of
non-limit ordinals less than $\delta$ which is cofinal in
$\delta$, i.e. $\sup\{ \eta_{\delta}(n) : n \in \omega \}=\delta$.
\item The ladder $\eta_{\delta}$ is a {\rm special ladder} if
there exists a sequence $0 < k^{\eta_{\delta}}_0 <
k^{\eta_{\delta}}_1 < \cdots < k^{\eta_{\delta}}_n < \cdots$ of
natural numbers such that
\begin{enumerate}
\item $\eta_{\delta}(k^{\eta_{\delta}}_n+i) + \omega =
\eta_{\delta}(k^{\eta_{\delta}}_n+j) + \omega$ for all $i,j <
k^{\eta_{\delta}}_{n+1} - k^{\eta_{\delta}}_n$; \item
$\eta_{\delta}(k^{\eta_{\delta}}_n) + \omega <
\eta_{\delta}(k^{\eta_{\delta}}_{n+1})$. \end{enumerate} for all
$n \in \omega$.
\end{enumerate}
\end{definition}

Note that the existence of $S$ and certain ladders on $S$ is
well-known. However, any limit ordinal $\delta$ of the form
$\delta=\alpha + \omega$ obviously does not allow the existence of
a special ladder. This is the reason why we have required that
$\omega^2$ divides $\delta$ for every $\delta \in S$, hence no
$\delta \in S$ can be of the form $\delta=\alpha + \omega$. For
$\delta \in S$ it is then easy to see how to obtain a
special ladder from a given ladder on $\delta$.\\

\begin{example}
\label{example1} The following are natural examples for a ladder
$\eta$ (on $\delta$) to be special:
\begin{enumerate}
\item Let $k_n^{\eta}=2n$ for all $n \in \omega$. Then $\eta$ is
special if and only if \[\eta(2n)+\omega = \eta(2n+1) + \omega <
\eta(2n+2);\] \item Let $k_n^{\eta}=n$ for all $n \in \omega$.
Then $\eta$ is special if and only if \[\eta(n)+\omega <
\eta(n+1).\]
\end{enumerate}
\end{example}

For $\delta \in S$ we let $\Delta_{\delta}$ be the set of all
special ladders on $\delta$ defined by
\[ \Delta_{\delta}= \{ \eta : \omega \rightarrow \delta
| \eta \text{ is a special ladder} \}. \]

We now collect ladder systems containing special ladders.

\begin{definition}
A system $\bar{\eta}= \left< \eta_{\delta} : \delta \in S \right>$
of (special) ladders is called a {\rm (special) ladder system on
$S$}.
\end{definition}

We put
\[ E = \{ \bar{\eta} : \bar{\eta} \text{ is
a special ladder system} \}.
\]
For later use we also define

\[ E_{\alpha} = \{ \bar{\eta}_{\alpha} : \bar{\eta} \in E, \bar{\eta}_{\alpha}=\left< \eta_{\delta} : \delta \in S \cap \alpha \right>\}
\]
for $\alpha < \omega_1$. On the set of special ladders we define
the {\it $\omega$-range} function as follows: \[ \rd :
\bigcup_{\delta \in S} \Delta_{\delta} \rightarrow \im (\rd), \eta
\mapsto \left< \eta(k^{\eta}_n) + \omega : n \in \omega \right>.
\]
Note that $\rd(\eta)$ determines all values of $\eta(n)+\omega$
$(n \in \omega)$ since the ladder is special. Moreover, if
$\bar{\eta} \in E$, then put $\rd(\bar{\eta})= \left<
\rd(\eta_{\delta}) : \delta \in S \right>$ and similarly
$\rd_{\alpha}(\bar{\eta})=
\rd_{\alpha}(\bar{\eta}_{\alpha})=\left<
\rd(\eta_{\delta}) : \delta \in S \cap \alpha \right>$ for $\alpha < \omega_1$.\\

Using the special ladder systems we can now define our desired
groups. Let $\bar{\eta} \in E$ be a special ladder system and put
$k_n^{\delta}=k_n^{\eta_{\delta}}$ for all $\delta \in S$ and $n
\in \omega$. Moreover, let $t_n^{\delta}=k_{n+1}^{\delta} -
k_n^{\delta}$ for all $n \in \omega$. We define a $\Q$-module
\[ F= \bigoplus\limits_{\beta < \omega_1} x_{\beta}\Q \oplus
\bigoplus\limits_{\delta \in S, n \in \omega} y_{\delta,n}\Q
\] freely generated (as a vectorspace) by the independent elements
$x_{\beta}$ $(\beta < \omega_1)$ and $y_{\delta,n}$ $(\delta \in
S, n \in \omega)$. Our desired group will be constructed as a
subgroup of $F$. Therefore, given a group $G \subseteq F$, we
define a {\it canonical $\aleph_1$-filtration} of $G$ by letting
\[ G^{\alpha} = \left< G \cap ( \{x_{\beta} : \beta < \alpha + \omega\}
\cup \{ y_{\delta,n} : \delta \in S \cap \alpha, n \in \omega \})
\right>_* \subseteq G \] for $\alpha < \omega_1$. Then $\{
G^{\alpha} : \alpha < \omega_1 \}$ is an increasing continuous
chain of pure subgroups of $G$ such that $G=\bigcup\limits_{\alpha
< \omega_1}G^{\alpha}$. For simplicity let
$y_{\delta}=y_{\delta,0}$ for $\delta \in S$. Let $\psi: \omega
\rightarrow \omega$ be a set function with $\psi(n) \not=0$ for
all $n \in \omega$ and choose integers $a_l^{\delta,n}$ for $l <
t_n^{\delta}$ such that $\gcd(a_l^{\delta,n} : l <
t_n^{\delta})=1$ for all $n \in \omega$ and $\delta \in S$. We
define elements $z_{\delta,n} \in F$ via
\begin{equation}
\label{definition} z_{\delta,n}=
\prod_{i=0}^{n-1}\frac{1}{\psi(i)}y_{\delta} +
\sum_{i=0}^{n-1}\prod_{j=i}^{n-1}\frac{1}{\psi(j)}\left(\sum_{l <
t_i^{\delta}}a_l^{\delta,i}x_{\eta_{\delta}(k_i^{\delta}+l)}
\right)
\end{equation} for $\delta \in S$ and $n
\in \omega$. Furthermore we put $z_{\delta,0}=y_{\delta}$ and let $\bar{a}=\left<a_l^{\delta,n} : l < t_n^{\delta}, n \in \omega, \delta \in S \right>$.\\

Let $G_{\bar{\eta}}^{\psi,\bar{a}}=\left< x_{\beta}, z_{\delta,n} : \beta <
\omega_1, \delta \in S, n \in \omega \right> \subseteq F$. Then
easy calculations show that the only relations satisfied by the
generators of $G_{\bar{\eta}}^{\psi,\bar{a}}$ are
\begin{equation}
\label{main}
 \psi(n)z_{\delta,n+1} = z_{\delta,n} + \sum_{l <
t_n^{\delta}}a_l^{\delta,n}x_{\eta_{\delta}(k_n^{\delta}+l)}
\end{equation}
for $\delta \in S$ and $n \in \omega$. To simplify notation we shall omit in the sequel the superscript $(\psi,\bar{a})$ since the function $\psi$ and the vector $\bar{a}$ of integers will always be clear from the context. However, the reader should keep in mind that for every ladder system $\bar{\eta}$ the group $G_{\bar{\eta}}=G_{\bar{\eta}}^{\psi,\bar{a}}$ always depends on the additional parameters $\psi$ and $\bar{a}$. We consider Example
\ref{example1} again.

\begin{example} The following hold:
\label{example14}
\begin{enumerate}
\item Let $\bar{\eta}$ be a special ladder system consisting of
special ladders as defined in Example \ref{example1} (i) and
choose $a_0^{\delta,n}=1$, $a_1^{\delta,n}=-1$ for all $\delta
\in S$ and $n \in \omega$. Then $G_{\bar{\eta}}$ satisfies the
following relations
\[ \psi(n)z_{\delta,n+1} = z_{\delta,n} + x_{\eta_{\delta}(2n)} -
x_{\eta_{\delta}(2n+1)}\] \item Let $\bar{\eta}$ be a special
ladder system consisting of special ladders as defined in Example
\ref{example1} (ii) and choose $a_0^{\delta,n}=1$ for all $\delta
\in S$ and $n \in \omega$. Then $G_{\bar{\eta}}$ satisfies the
following relations
\[ \psi(n)z_{\delta,n+1} = z_{\delta,n} + x_{\eta_{\delta}(n)}\]
\end{enumerate}
\end{example}

We now prove some properties of the constructed groups
$G_{\bar{\eta}}$.

\begin{lemma}
Let $\bar{\eta} \in E$. Then the group $G_{\bar{\eta}}$ is a
torsion-free $\aleph_1$-separable abelian group of size $\aleph_1$
with $\Gamma(G_{\bar{\eta}})=\tilde{S}$.
\end{lemma}

\proof Let $\bar{\eta} \in E$ be a special ladder system. Clearly
the group $G_{\bar{\eta}}$ is a torsion-free group of cardinality
$\aleph_1$. We first prove that $G_{\bar{\eta}}$ is
$\aleph_1$-free. Therefore, let $H$ be a finite rank subgroup of
$G_{\bar{\eta}}$. Then there exists a finite subset
\[
T \leq \{x_{\beta}: \beta < \omega_1 \} \cup \{ z_{\delta,n} :
\delta \in S, n \in \omega \} \] such that
\[
H \subseteq \left< t : t \in T \right>_* \subseteq G_{\bar{\eta}}.
\]
Let $T_S=\{ \delta \in S : z_{\delta,n} \in T \text{ for some } n \in \omega \}$.
By enlarging $T$ we may assume that there exists an integer $m$
such that
\begin{itemize}
\item for $\delta \in T_S$ we have $z_{\delta,n} \in T$ if and only if $n \leq m$; \item for
$y_{\delta}=z_{\delta,0} \in T$ we have $x_{\eta_{\delta}(n)} \in T$ if and
only if $n < k_{m+1}^{\delta}$.
\end{itemize} Then using equation (\ref{main}) it is not hard to
see that $\left< t : t \in T \right>_*$ is freely generated by the
elements $\{ z_{\delta,m} : y_{\delta} \in T \} \cup \{
x_{\eta_{\delta}(n)} : n < k_{m+1}^{\delta}, y_{\delta} \in T \}$.
Thus $H$
is free and therefore $G_{\bar{\eta}}$ is $\aleph_1$-free.\\
It remains to prove that $G_{\bar{\eta}}$ is $\aleph_1$-separable. Therefore
let $\{ G_{\bar{\eta}}^{\alpha} : \alpha < \omega_1 \}$ be the
canonical $\aleph_1$-filtration of $G_{\bar{\eta}}$. We shall now
define for all $\nu \not\in S$ a projection $\pi_{\nu} :
G_{\bar{\eta}} \rightarrow G_{\bar{\eta}}^{\nu}$ such that
$\pi_{\nu} \restriction_{G_{\bar{\eta}}^{\nu}} = id
\restriction_{G_{\bar{\eta}}^{\nu}}$. Let $\nu \not\in S$ be
given. For every $\mu \geq \nu + \omega$ let
$\pi_{\nu}(x_{\mu})=0$; for $\delta \in S$ with $\delta > \nu$ let
$n_{\delta}$ be maximal with
$\eta_{\delta}(k_{n_{\delta}}^{\delta}) < \nu$. Hence
$\eta_{\delta}(k^{\delta}_{n_{\delta}}+i) < \nu + \omega$ for all
$i < t_{n_{\delta}}^{\delta}$ and
$\eta_{\delta}(k_{n_{\delta}+1}^{\delta}) > \nu + \omega$. Let
$\pi_{\nu}(z_{\delta,n})=0$ for all $n \geq n_{\delta}$. Moreover,
put
\[ \pi_{\nu}(y_{\delta})= -
\sum_{i=0}^{n_{\delta}-1}\prod_{j=0}^{i-1} \psi(j+1)\sum_{l <
t_i^{\delta}}a_l^{\delta,i}x_{\eta_{\delta}(k_i^{\delta}+l)}
\] and finally
\[
\pi_{\nu}(z_{\delta,n})=-\sum_{i=n}^{n_{\delta}-1}\prod_{j=n}^{i-1}\psi(j+1)\sum_{l
< t_i^{\delta}}a_l^{\delta,i}x_{\eta_{\delta}(k_i^{\delta}+l)}
\] for all $n < n_{\delta}$. Letting $\pi_{\nu}
\restriction_{G_{\bar{\eta}}^{\nu}} = id
\restriction_{G_{\bar{\eta}}^{\nu}}$ it is now straightforward to
check that $\pi_{\nu}$ is a well-defined homomorphism as claimed
using equation (\ref{main}). Finally,
$\Gamma(G_{\bar{\eta}})=\tilde{S}$ follows immediately. \qed\\

We now prove that a special ladder system is enough {\it
separated}.

\begin{lemma}
\label{integers} Let $\bar{\eta} \in E$ and $\alpha < \omega_1$.
Then there exists a sequence of integers $\left< m_{\delta} :
\delta \in S \cap \alpha \right>$ such that the sets $\{
\eta_{\delta}(k_n^{\delta}) + \omega : n \geq m_{\delta} \}$
($\delta \in S \cap \alpha$) are pairwise disjoint. In
particular, the sets $\{ \eta_{\delta}(k_n^{\delta}+i) : n \geq
m_{\delta}, i < t_n^{\delta} \}$ ($\delta \in S \cap \alpha$) are
pairwise disjoint.
\end{lemma}

\proof Let $\bar{\eta} \in E$ and $\alpha < \omega_1$ be given.
Since $\alpha$ is countable we may enumerate $S \cap \alpha$ by
$\omega$, say $S \cap \alpha=\{ \delta_k : k \in \omega \}$. We
shall now define inductively the sequence $\left< m_{\delta_k} : k
\in \omega \right>$ such that for every $k \in \omega$ the sets
\begin{equation}
\label{disjoint} \{\eta_{\delta_j}(k_n^{\delta_j}) + \omega : n
\geq m_{\delta_j} \} \quad (j \leq k) \text{ are pairwise
disjoint.}
\end{equation} We start with $k=1$, hence $\eta_{\delta_0}$ and
$\eta_{\delta_1}$ are given. In order to carry on the induction we shall prove a stronger
result. Let $m_{\delta_0}$ be fixed but arbitrary. We claim that there is $m_{\delta_1}$ such that (\ref{disjoint}) holds for $k=1$. Assume first that $\delta_0 < \delta_1$. Since $S \subseteq
\lim(\omega_1)$ and $\omega^2 | \delta$ for all $\delta \in S$ we obtain $\delta_1 > \delta_0 + \omega$. Hence it is easy
to see that $m_{\delta_1}$ exists such that
(\ref{disjoint}) is satisfied for $k=1$ because $\eta_{\delta_1}$
is a ladder
with $\sup(\im(\eta_{\delta_1})) =\delta_1$.\\
If $\delta_1 < \delta_0$, then $\delta_0 > \delta_1 + \omega$. Thus there is $m'_{\delta_1}$ such that
$\{\eta_{\delta_0}(k_n^{\delta_0}) + \omega : n \geq m'_{\delta_1} \}$ and $\{
\eta_{\delta_1}(k_n^{\delta_1}) + \omega : n \geq m'_{\delta_1} \}$ are disjoint. Increasing $m'_{\delta_1}$ sufficiently we obtain $m_{\delta_1} \geq m'_{\delta_1}$ such that (\ref{disjoint}) holds.\\

The inductive step is now immediate. Given $k$ such that
$m_{\delta_0}, m_{\delta_1}, \cdots, m_{\delta_{k-1}}$ satisfy
(\ref{disjoint}) we obtain integers
$s_j$ for $j < k$ such that $\{ \eta_{\delta_j}(k_n^{\delta_j}) +
\omega : n \geq m_{\delta_j}\}$ and $\{ \eta_{\delta_k}(k_n^{\delta_k}) +
\omega : n \geq s_j \}$ are pairwise disjoint for every $j < k$.
Choosing $m_{\delta_k}=\max\{ s_j : j < k \}$ we satisfy
(\ref{disjoint}).
\qed\\

Note that Lemma \ref{integers} gives the same sequence of integers
for different $\bar{\eta}, \bar{\nu} \in E$ if
$\rd(\bar{\eta})=\rd(\bar{\nu})$. Nevertheless, the next lemma
shows that special ladder systems $\bar{\eta}, \bar{\nu} \in E$
with $\rd(\bar{\eta})=\rd(\bar{\nu})$ do not overlap very much.

\begin{lemma}
\label{overlap} Let $\bar{\eta}, \bar{\nu} \in E$ and $\alpha <
\omega_1$ such that $\rd(\bar{\eta})=\rd(\bar{\nu})$. Moreover,
let $\left< m_{\delta} : \delta \in S \cap \alpha \right>$ be the
sequence from Lemma \ref{integers}. If
$\eta_{\delta}(k_n^{\eta_{\delta}}+j)=\nu_{\delta'}(k_m^{\nu_{\delta'}}+i)$
for some $n \geq m_{\delta}$, $m \geq m_{\delta'}$, and $i <
t_m^{\nu_{\delta'}}$, $j <
t_n^{\nu_{\delta}}$. Then $\delta=\delta'$ and $m=n$.
\end{lemma}

\proof Assume that
$\eta_{\delta}(k_n^{\eta_{\delta}}+j)=\nu_{\delta'}(k_m^{\nu_{\delta'}}+i)$
for some $n \geq m_{\delta}$, $m \geq m_{\delta'}$, and $i <
t_m^{\nu_{\delta'}}$, $j <
t_n^{\nu_{\delta}}$. Then
\[ \eta_{\delta}(k_n^{\eta_{\delta}}+j) + \omega=  \eta_{\delta}(k_n^{\eta_{\delta}}) + \omega = \nu_{\delta'}(k_m^{\nu_{\delta'}}+i) + \omega
=\nu_{\delta'}(k_m^{\nu_{\delta'}}) + \omega =
\eta_{\delta'}(k_m^{\eta_{\delta'}}) + \omega \] since
$\rd(\bar{\eta})=\rd(\bar{\nu})$. Thus $\delta=\delta'$ by Lemma
\ref{integers}. Moreover, $m=n$ follows since $\eta_{\delta}$ is
a special ladder. \qed\\

Recall that two groups $A$ and $B$ of
cardinality $\aleph_1$ are called {\it filtration-equivalent} if
they have filtrations $\{ A_{\nu} : \nu \in \omega_1 \}$ and $\{
B_{\nu} : \nu \in \omega_1 \}$ respectively such that for all $\nu
\in \omega_1$, there is an isomorphism $\Phi_{\nu}: A_{\nu}
\rightarrow B_{\nu}$ satisfying $\Phi_{\nu}[A_{\mu}]=B_{\mu}$ for all
$\mu \leq \nu$. Such an isomorphism is called {\it
level-preserving}. Note that we do not require that $\Phi_{\tau}$
extends $\Phi_{\nu}$ when $\tau \geq \nu$ and that
filtration-equivalent groups $A$ and $B$ are also
quotient-equivalent, i.e. for all $\nu \in \omega_1$ we have
$A_{\nu +1}/ A_{\nu} \cong B_{\nu + 1}/B_{\nu}$.

\begin{proposition}
\label{mainprop} Let $\bar{\eta}, \bar{\nu} \in E$ such that
$\rd(\bar{\eta})=\rd(\bar{\nu})$. Then the groups $G_{\bar{\eta}}$
and $G_{\bar{\nu}}$ are filtration-equivalent.
\end{proposition}

\proof Let $\bar{\eta}$ and $\bar{\nu}$ be given. By construction
we have
\[ G_{\bar{\eta}}=\left< x_{\beta}, z_{\delta,n} : \beta <
\omega_1, \delta \in S, n \in \omega \right> \] and
\[ G_{\bar{\nu}}=\left< x_{\beta}, w_{\delta,n} : \beta <
\omega_1, \delta \in S, n \in \omega \right> \] such that the
elements $z_{\delta,n}$ and $w_{\delta,n}$ $(\delta \in S, n \in
\omega)$ are defined as in (\ref{definition}) for $\bar{\eta}$ and
$\bar{\nu}$ respectively. Hence, the only relations satisfied in
$G_{\bar{\eta}}$ and $G_{\bar{\nu}}$ are the relations in equation
(\ref{main}). Since filtration-equivalence is a transitive
property it suffices to assume that $G_{\bar{\eta}}$ is of the
simplest form $(k_n^{\eta_{\delta}}=n, t_n^{\eta_{\delta}}=1,
a_0^{\eta_{\delta},n}=1$), hence
\[\psi(n)z_{\delta,n+1} = z_{\delta,n} +
x_{\eta_{\delta}(n)}\] and
\[ \psi(n)w_{\delta,n+1} = w_{\delta,n} + \sum_{l <
t_n^{\delta}}a_l^{\delta,n}x_{\nu_{\delta}(k_n^{\delta}+l)}. \]
Note that the parameters $k_n^{\delta}$, $t_n^{\delta}$ and
$a_l^{\delta,n}$ depend on $\nu_{\delta}$. Moreover, we shall
assume for simplicity and without loss of generality that
$a_0^{\delta,n}=1$ for every $\delta \in S$, $n \in \omega$ since $\gcd(a_l^{\delta,n} : l < t_n^{\delta})=1$.
Hence we may replace the basis element
$x_{\nu_{\delta}(k_n^{\delta})}$ by the new basis element
$\sum_{l <
t_n^{\delta}}a_l^{\delta,n}x_{\nu_{\delta}(k_n^{\delta}+l)}.$

Let $G_{\bar{\eta}}=\bigcup\limits_{\alpha <
\omega_1}G^{\alpha}_{\bar{\eta}}$ and
$G_{\bar{\nu}}=\bigcup\limits_{\alpha <
\omega_1}G^{\alpha}_{\bar{\nu}}$ be the canonical
$\aleph_1$-filtrations of $G_{\bar{\eta}}$ and $G_{\bar{\nu}}$
respectively. For each $\alpha < \omega_1$ we now define a
level-preserving isomorphism from $G^{\alpha}_{\bar{\eta}}$ onto
$G^{\alpha}_{\bar{\nu}}$. Let $\alpha < \omega_1$ be fixed. Since
by assumption $\rd(\bar{\eta})=\rd(\bar{\nu})$ we may choose a
sequence $\left< m_{\delta} : \delta \in S \cap \alpha \right>$ as
in Lemma \ref{integers} for $\bar{\eta}$ and $\bar{\nu}$
simultaneously.
Let
$\hat{\pi}:G^{\alpha}_{\bar{\eta}} \rightarrow
G^{\alpha}_{\bar{\nu}}$ be defined via
\begin{itemize}
\item $\hat{\pi}(x_{\eta_{\delta}(n)}) = \sum_{l <
t_n^{\delta}}a_l^{\delta,n}x_{\nu_{\delta}(k_n^{\delta}+l)}
\textit{ for all } n \geq m_{\delta}, \delta \in S \cap \alpha$

\item $\hat{\pi}(x_{\nu_{\delta}(k_n^{\delta})}) = x_{\eta_{\delta}(n)} \textit{ for all } n
\geq m_{\delta}, \delta \in S \cap \alpha \text{ if }
\eta_{\delta}(n) \not= \nu_{\delta}(k_n^{\delta})$

\item $\hat{\pi}(x_{\beta}) = x_{\beta} \textit{ for every } \beta < \alpha + \omega \textit{ otherwise and}
$
\item $\hat{\pi}(z_{\delta,n}) = w_{\delta,n} \textit{ for all } n
\geq m_{\delta}, \delta \in S \cap \alpha. $ \end{itemize}

Recursively we may define $\hat{\pi}(z_{\delta,n})$ for $n <
m_{\delta}$ $(\delta \in S \cap \alpha$) using the definition of
$\hat{\pi}$ on $x_{\beta}$ $(\beta < \alpha + \omega)$ and on
$z_{\delta,m_{\delta}}$. By the choice of the sequence $\left<
m_{\delta} : \delta \in S \cap \alpha \right>$ it is now easy to
see that $\hat{\pi}$ is a level preserving isomorphism from
$G^{\alpha}_{\bar{\eta}}$ onto $G^{\alpha}_{\bar{\nu}}$ and hence
the groups $G_{\bar{\eta}}$ and $G_{\bar{\nu}}$ are filtration-equivalent. \qed

\section{The Consistency Result}
\relax From now on we let $\psi: \omega \rightarrow \omega$ be given by
$\psi(n)=n!$ with the convention that $0!=1$. In order to force that the group $G_{\bar{\eta}}$ is
a Whitehead group we recall the definition of the uniformization
property.

\begin{definition}
If $\lambda$ is a cardinal and $\bar{\eta}$ is a ladder system on
$S$ we say that $\bar{\eta}$ has {\it $\lambda$--uniformization}
if for every family $\{ c_{\delta} : \delta \in S\}$ of colors
$c_{\delta}: \omega \rightarrow \lambda$, there exist $\Psi:
\omega_1 \rightarrow \lambda$ and $\Psi^*: S \rightarrow \omega$
such that $\Psi(\eta_{\delta}(n))=c_{\delta}(n)$ for all
$\Psi^*(\delta)\leq n$ and $\delta \in S$.
\end{definition}

The following lemma is by now standard (compare \cite[Chapter
XIII, Proposition 0.2]{EkMe}). However, the construction in
\cite[Chapter XIII, Section 0]{EkMe} is slightly different from
our construction since $x_{\eta_{\delta}(k_n^{\eta_{\delta}}+i)}$
($i < t_n^{\eta_{\delta}})$ appear in equation (\ref{main}) at the
same time. Therefore, we give the adjusted proof of the next
lemma in a particular case for the convenience of the reader. However, we shall only apply it for $\bar{\eta}$ of the simplest form as in Example \ref{example14} (ii).

\begin{lemma}
\label{uniform1} If $\bar{\eta}$ is a ladder system which has
$\aleph_0$-uniformization, then the group $G_{\bar{\eta}}$
satisfies $\Ext(G_{\bar{\eta}},N)=0$ for every countable abelian
group $N$. If $\bar{\eta}$ has $2$-uniformization then
$G_{\bar{\eta}}$ is a Whitehead group.
\end{lemma}

\proof Let $N$ be a countable abelian group. For simplicity we
shall assume the setting of Example \ref{example14} (i). The general proof is similiar. By
construction we may regard $G_{\bar{\eta}}$ as the quotient $P/Q$
of the free group $P=\bigoplus\limits_{\beta <
\aleph_1}x_{\beta}\Z \oplus \bigoplus\limits_{\delta \in S, n \in
\omega}z_{\delta,n}\Z$ and its subgroup $Q$ generated by the
elements \[ g_{\delta,n}=n!z_{\delta,n+1} - z_{\delta,n} -
x_{\eta_{\delta}(2n+1)} + x_{\eta_{\delta}(2n)} \] for $\delta \in
S$ and $n \in \omega$. In order to show that
$\Ext(G_{\bar{\eta}},N)=0$ it therefore suffices to prove that
every homomorphism $\varphi: Q \rightarrow N$ has an extension
$\tilde{\varphi} : P \rightarrow N$. Thus let $\varphi \in
\Hom(G_{\bar{\eta}},N)$ be given. We fix a bijection $b: N
\rightarrow \aleph_0$ and define $c_{\delta}: \omega \rightarrow
\omega$ for $\delta \in S$ as follows: Let $n \in \omega$ and put
\[
c_{\delta}(2n)=b(\varphi(g_{\delta,n})) \text{ and }
c_{\delta}(2n+1)=b(2\varphi(g_{\delta,n})).
\]
By the uniformization property there exists $f: \omega_1 \rightarrow
\omega$ such that for all $\delta \in S$ there exists $k_{\delta}
\in \omega$ such that
\[
c_{\delta}(n)=f(\eta_{\delta}(n)) \text{ for all } n > k_{\delta}.
\]
We define $\tilde{\varphi} : P \rightarrow N$ as follows: Let
$\alpha < \omega_1$
\begin{itemize}
\item If $\alpha = \eta_{\delta}(n)$ for some $\delta \in S$ and
$n> k_{\delta}$ then put
$\tilde{\varphi}(x_{\alpha})=b^{-1}(f(\alpha))$; note that $\alpha
\not\in S$; \item If $\alpha \not\in S$ and $\alpha \not=
\eta_{\delta}(n)$ for any $\delta \in S$ and $n> k_{\delta}$ then
put $\tilde{\varphi}(x_{\alpha})=0$; \item if $\alpha \in S$ and
$2n > k_{\alpha}$ then put $\tilde{\varphi}(z_{\alpha,n})=0$;
\item if $\alpha \in S$ and $2n \leq k_{\alpha}$ then we define
$\tilde{\varphi}(z_{\alpha,n})$ inductively and distinguish the
following four cases:
\begin{itemize}
\item if $\eta_{\alpha}(2n)=\eta_{\delta}(k)$ for some $k >
k_{\delta}$ and $\eta_{\alpha}(2n+1)=\eta_{\delta'}(k')$ for some
$k'
> k_{\delta'}$ then put
\[
\tilde{\varphi}(z_{\alpha,n})=b^{-1}(f(\eta_{\alpha}(2n+1))) -
b^{-1}(f(\eta_{\alpha}(2n))) - \varphi(g_{\alpha,n}) +
n!\tilde{\varphi}(z_{\alpha,n+1});\] \item if
$\eta_{\alpha}(2n)=\eta_{\delta}(k)$ for some $k > k_{\delta}$ but
$\eta_{\alpha}(2n+1) \not=\eta_{\delta'}(k')$ for all $k'
> k_{\delta'}$ and $\delta' \in S$ then put
\[\tilde{\varphi}(z_{\alpha,n})= - b^{-1}(f(\eta_{\alpha}(2n))) -
\varphi(g_{\alpha,n}) + n!\tilde{\varphi}(z_{\alpha,n+1});\] \item
if $\eta_{\alpha}(2n) \not=\eta_{\delta}(k)$ for all $k >
k_{\delta}$ and $\delta \in S$ but
$\eta_{\alpha}(2n+1)=\eta_{\delta'}(k')$ for some $k'
> k_{\delta'}$ then put
\[\tilde{\varphi}(z_{\alpha,n})=b^{-1}(f(\eta_{\alpha}(2n+1))) -
\varphi(g_{\alpha,n}) + n!\tilde{\varphi}(z_{\alpha,n+1});\] \item
if $\eta_{\alpha}(2n) \not=\eta_{\delta}(k)$ for all $k >
k_{\delta}$ and $\delta \in S$ and also
$\eta_{\alpha}(2n+1)\not=\eta_{\delta'}(k')$ for all $k'
> k_{\delta'}$ and $\delta' \in S$ then put
\[\tilde{\varphi}(z_{\alpha,n})= - \varphi(g_{\alpha,n}) +
n!\tilde{\varphi}(z_{\alpha,n+1}).\]
\end{itemize}
\end{itemize}
It now remains to show that $\tilde{\varphi}$ is an extension of
$\varphi$, hence satisfies
$\tilde{\varphi}(g_{\alpha,n})=\varphi(g_{\alpha,n})$ for all
$\alpha \in S$ and $n \in \omega$. Clearly we have
\[
\tilde{\varphi}(g_{\alpha,n})=n!\tilde{\varphi}(z_{\alpha,n+1}) -
\tilde{\varphi}(z_{\alpha,n})-\tilde{\varphi}(x_{\eta_{\alpha}(2n)})
+ \tilde{\varphi}(x_{\eta_{\alpha}(2n+1)}). \] If $\alpha \in S$
and $2n
> k_{\alpha}$ then
\[\tilde{\varphi}(x_{\eta_{\alpha}(2n)})=b^{-1}(f(\eta_{\alpha}(2n)))=b^{-1}(c_{\alpha}(2n))=\varphi(g_{\alpha,n})\]
and similarly
$\tilde{\varphi}(x_{\eta_{\alpha}(2n+1)})=2\varphi(g_{\alpha,n})$.
Furthermore,
$\tilde{\varphi}(z_{\alpha,n})=\tilde{\varphi}(z_{\alpha,n+1})=0$
and hence \[
\tilde{\varphi}(g_{\alpha,n})=-\varphi(g_{\alpha,n})+2
\varphi(g_{\alpha,n})=\varphi(g_{\alpha,n}). \] All other cases
can be checked similarly by easy calculations and are therefore
left to the reader.\\
The second statement follows similarly using \cite[Chapter XIII,
Lemma 0.7]{EkMe}
 \qed\\

Similarly, we can prove the next result which is essentially
\cite[Chapter XIII, Proposition 0.6]{EkMe}. Recall that a ladder system $\bar{\eta}$ is called {\it tree-like} if $\eta_{\delta}(n)=\eta_{\delta'}(m)$ for some $\delta, \delta' \in
S$ and $n,m \in \omega$ implies $m=n$ and $\eta_{\delta}(k)=\eta_{\delta'}(k)$ for all $k \leq n$.

\begin{lemma}
\label{uniform2} Let $\bar{\eta}$ be a special tree-like ladder
system. If $G_{\bar{\eta}}$ satisfies
$\Ext(G_{\bar{\eta}},\Z^{(\omega)})=0$, then $\bar{\eta}$ has
$\aleph_0$-uniformization. In particular, if $G_{\bar{\eta}}$ is a
Whitehead group, then $\bar{\eta}$ has $2$-uniformization.
\end{lemma}

\proof As in the proof of Lemma \ref{uniform1} we shall assume
for simplicity the setting of Example \ref{example14} and let
$G_{\bar{\eta}}=P/Q$. Let $\{ a_{nmj} : n,m,j \in \omega \}$ be a
basis of $\Z^{(\omega)}$. Given an $\omega$-coloring $\{
c_{\delta} : \delta \in S \}$ define $\varphi: Q \rightarrow
\Z^{(\omega)}$ by
\[ \varphi(g_{\delta,n})=a_{nc_{\delta}(2n+1)c_{\delta}(2n+2)}.\]
By hypothesis there exists an extension of $\varphi$ to
$\tilde{\varphi}:P \rightarrow \Z^{(\omega)}$. Define
$\Psi^*(\delta)$ to be the least integer $n' > 4$ such that
\[ \varphi(z_{\delta,0}) \in \left< a_{lmj} : l < n',m,j \in \omega
\right>.\] It suffices to show that if
$\eta_{\delta}(k)=\eta_{\gamma}(k)$ where $k \geq
\Psi^*(\delta),\Psi^*(\gamma)$ then $c_{\delta}(k)=c_{\gamma}(k)$.
In this case $\Psi(\eta_{\delta}(k)))=c_{\delta}(k)$ when $k \geq
\Psi^*(\delta)$ is as required. Thus let $k \geq \Psi^*(\delta),
\Psi^*(\gamma)$ and $\eta_{\delta}(k)=\eta_{\gamma}(k)$. Then
$k=2n+1$ or $2n+2$ for some $2 \leq n \in \omega$. Let
$\tilde{\varphi}'$ be the composition of $\tilde{\varphi}$ with
the projection of $\Z^{(\omega)}$ onto $\left< a_{kmj} : m,j \in
\omega \right>$. Then
\[\tilde{\varphi}'(z_{\delta,0})=\tilde{\varphi}'(z_{\gamma,0})=0.\]
Since $\bar{\eta}$ is tree-like we have
$x_{\eta_{\delta}(s)}=x_{\eta_{\gamma}(s)}$ for all $s \leq k$.
Using this and the fact that
$\tilde{\varphi}'(g_{\delta,n'})=\tilde{\varphi}'(g_{\gamma,n'})=0$
for all $n' < k$ we can show by induction that
\[\tilde{\varphi}'(z_{\delta,n})=\tilde{\varphi}'(z_{\gamma,n}).\]
Hence
\[ a_{nc_{\delta}(2n+1)c_{\delta}(2n+2)} -
a_{nc_{\gamma}(2n+1)c_{\gamma}(2n+2)}= \varphi(g_{\delta,n}) -
\varphi(g_{\gamma,n}) =
n!(\tilde{\varphi}'(z_{\gamma,n+1})-\tilde{\varphi}'(z_{\delta,n+1}).
\]
Therefore $n!$ divides $a_{nc_{\delta}(2n+1)c_{\delta}(2n+2)} -
a_{nc_{\gamma}(2n+1)c_{\gamma}(2n+2)}$; so
$a_{nc_{\delta}(2n+1)c_{\delta}(2n+2)}$ must equal $
a_{nc_{\gamma}(2n+1)c_{\gamma}(2n+2)}$ since they are basis
elements and hence $c_{\delta}(k)=c_{\gamma}(k)$ since either
$k=2n+1$ or $k=2n+2$.\\
The second statement follows similarly with the appropriate
adjustments and \cite[Chapter XIII, Proposition 0.6]{EkMe}. \qed\\

We are now ready to prove the main theorem. Therefore let
$\bar{\nu}$ be a special ladder system and $\bar{c}=\left<
c_{\delta} : \omega \rightarrow \{0,1\} | \delta \in S \right>$ a
sequence of colors. We define a group $H_{\bar{\nu},\bar{c}}$ as
follows. Similar to the group $G_{\bar{\nu}}$ constructed in the
previous section we let $F'$ be the free abelian group
\[ F'= \hat{w}\Z \oplus \bigoplus\limits_{\alpha < \omega_1}\hat{x}_{\alpha} \oplus \bigoplus\limits_{\delta \in S, n \in \omega} y'_{\delta,n}\Z \]
and $H_{\bar{\nu},\bar{c}}$ be the subgroup of $F'$ generated by
\[ H_{\bar{\nu},\bar{c}}= \left< \hat{x}_{\beta}, \hat{z}_{\delta,n}: \beta < \omega_1, \delta \in S, n \in \omega \right> \subseteq_* F', \]
where the $\hat{z}_{\delta,n}$ are chosen subject to the relations
\[n!\hat{z}_{\delta,n+1} = \hat{z}_{\delta,n} + \sum\limits_{i < t_n^{\delta}} a_l^{\delta,n} \hat{x}_{\nu_{\delta}(k_n^{\delta} + i)} + c_{\delta}(n)\hat{w} \]
for $\delta \in S$ and $n \in \omega$. We define a natural mapping
$h_{\bar─{\nu},\bar{c}}: H_{\bar{\nu},\bar{c}} \rightarrow
G_{\bar{\nu}}$ via
\begin{itemize}
\item $h_{\bar─{\nu},\bar{c}}(\hat{x}_{\beta})=x_{\beta}$ for all $\beta < \omega_1$;
\item $h_{\bar─{\nu},\bar{c}}(\hat{z}_{\delta,n})=z_{\delta,n}$ for all $\delta \in S$ and $n \in \omega$;
\item $h_{\bar─{\nu},\bar{c}}(\hat{w})=0$.
\end{itemize}
Obviously, the kernel of $h_{\bar─{\nu},\bar{c}}$ is isomorphic
to $\Z$, in fact $\ker(h_{\bar─{\nu},\bar{c}})=\hat{w}\Z$. Thus
$h_{\bar─{\nu},\bar{c}}$ induces a short exact sequence
\begin{equation}
\tag{E}
\label{E}
0 \longrightarrow \Z \longrightarrow H_{\bar{\nu},\bar{c}} \longrightarrow G_{\bar{\nu}} \longrightarrow 0.
\end{equation}
As for $G_{\bar{\nu}}$ we also define a filtration for $H_{\bar{\nu},\bar{c}}$ by letting
\[ H_{\bar{\nu},\bar{c}}^{\alpha}=\left< H_{\bar{\nu},\bar{c}} \cap (\{z, \hat{x}_{\beta} : \beta < \alpha + \omega \} \cup \{ y'_{\delta,n} : \delta \in S \cap \alpha, n \in \omega \}) \right>_* \subseteq H_{\bar{\nu},\bar{c}} \]
for $\alpha < \omega_1$.\\
The idea for proving the main theorem is to build an extension
model of ZFC in which GCH holds and in which we can force two
special ladder systems $\bar{\eta}$ and $\bar{\nu}$ with
$\rd(\bar{\eta})=\rd(\bar{\nu})$ such that $\bar{\eta}$ has the
$2$-uniformization property, hence $G_{\bar{\eta}}$ is a Whitehead
group but at the same time we force a coloring $\bar{c}$ such that
the sequence $(\ref{E})$ does not split, hence $G_{\bar{\nu}}$ is
not a Whitehead group. For notational reasons we call a special
ladder system of the simplest form as in Example \ref{example14}
(ii) a {\it simple special ladder system}.

\begin{theorem}
\label{maintheorem} There exists a model of $ZFC$ in which $GCH$
holds and for some special ladder systems $\bar{\eta},
\bar{\nu} \in E$ with $\rd(\bar{\eta})=\rd(\bar{\nu})$, the group
$G_{\bar{\eta}}$ is a Whitehead group but $G_{\bar{\nu}}$ is not.
\end{theorem}

\proof Essentially the proof is given in \cite{Sh2} (see also
\cite{Sh1} and \cite{Sh3}). Therefore we only recall the basic
steps of the proof. Suppose we start with a ground model $V$ in
which $GCH$ holds. Let $\bar{\eta}$ be a (special) ladder system.
It was shown in \cite[Theorem 1.1]{Sh1} that there exists a
forcing notion $(P, \leq)$ such that:
\begin{itemize}
\item $|P|=\aleph_2$, $P$ satisfies the $\aleph_2$-chain condition
and adds no new sequences of length $\omega$; hence, if $V$
satisfies $GCH$, then also the extension model $V^P$ satisfies
$GCH$; \item every stationary set remains stationary in $V^P$;
\item $\bar{\eta}$ has the $2$-uniformization property (even the
$\aleph_0$-uniformization property (see \cite[Theorem 2.1]{Sh1})).
\end{itemize}
The forcing notion $(P,\leq)$ was obtained by a countable support
iteration (of length $\aleph_2$); at each step using a basic
forcing extension and taking inverse limits at stages of
cofinality $\aleph_0$. We briefly recall the basic iteration. Let
$\bar{c}=\left< c_{\delta} : \delta \in S \right>$ be a system of
colors which has to be uniformized. Here each $c_{\delta}: \omega
\rightarrow 2$. Define $P_{\bar{c}}$ as the set of all functions
$f$ such that
\begin{enumerate}
\item $f: \alpha \rightarrow 2$ for some $\alpha < \omega_1$;
\item for all $\delta \leq \alpha, \delta \in S$ there is
$n_{\delta}$ such that $f(\eta_{\delta}(m))=c_{\delta}(m)$ for all
 $m \geq n_{\delta}$.
 \end{enumerate}
$P_{\bar{c}}$ is ordered naturally and it is easy to see that the
set $D_{\alpha}=\{ f \in P_{\bar{c}} : \alpha \subseteq \dom(f)
\}$ is dense for every $\alpha < \omega_1$ and hence a generic
filter will give the desired unifying $\Psi^*$.\\

Now, assume that $V \vDash GCH$ is given. We shall define a
countable support iteration $\bar{Q}=\left< P_{\alpha}, \overset{\cdot}{Q}_{\alpha}
: \alpha < \omega_2 \right>$ as follows: We start with an initial
forcing (compare also \cite{EkSh})
\begin{definition}
Let $P_0$ consist of all triples $\left< \bar{\eta}, \bar{\nu},
\bar{c} \right>$ such that for some $\alpha < \omega_1$ we have
 \begin{itemize} \item $\bar{\eta},
\bar{\nu} \in E_{\alpha}$ are special ladder systems on $S \cap
\alpha$ \item  $\bar{\eta}$ is simple \item $\bar{c}=\left<
c_{\delta}: \omega \rightarrow 2 | \delta \in S \cap \alpha
\right>$ \item
$\rd_{\alpha}(\bar{\eta})=\rd_{\alpha}(\bar{\nu})$.
\end{itemize}
\end{definition}
We may think of the conditions in $P_0$ as partial special ladder
systems on $S \cap \alpha$ for some $\alpha < \omega_1$ and a
corresponding partial coloring. It is easy to check that a
$P_0$-generic filter $G$ gives a pair of special ladder systems
$\bar{\eta}, \bar{\nu}$ on $S$ (in the extension model $V[G]$)
with the same $\omega$-range and a global coloring $\bar{c}$.
Moreover, $\bar{\eta}$ will be simple. Let $\tilde{\bar{\eta}},
\tilde{\bar{\nu}}$ and $\tilde{\bar{c}}$, $\tilde{c}_{\delta}$ be
the corresponding $P_0$-names which are defined naturally. Note
that $P_0$ is $\omega$-closed and satisfies the $\aleph_2$-chain
condition, so GCH holds in $V[G]$ since it holds in $V$. Applying
the forcing described above to $V'=V[G]$ we can force that in
$(V')^P$ the ladder system $\tilde{\bar{\eta}}$ has the
$2$-uniformization property and hence the group
$G_{\tilde{\bar{\eta}}}$ is a Whitehead group. Here, we let
$P=\left< P_{\alpha}, \overset{\cdot}{Q}_{\alpha} : 1 \leq \alpha < \omega_2
\right>$, hence in $V^{\bar{Q}}=(V')^P$ the generalized continuum
hypothesis holds. We have to show that the group
$G_{\tilde{\bar{\nu}}}$ is not a Whitehead group. As indicated
this shall be done by showing that the sequence $(\ref{E})$ cannot
be forced to split.\\
For the sake of contradiction assume that $(\ref{E})$ splits.
Hence for some $p \in \bar{Q}$ we have
\[ p \Vdash " f_{\tilde{\bar{\nu}}, \tilde{\bar{c}}} \in \Hom(G_{\tilde{\bar{\nu}}},H_{\tilde{\bar{\nu}}, \tilde{\bar{c}}})  \text{ is an inverse of } h_{\tilde{\bar{\nu}},\tilde{\bar{c}}}. "
\] Since $\bar{Q}$ satisfies the $\aleph_2$ chain condition we can
replace
$\bar{Q}$ by $P_{\alpha}$ for some $\alpha < \omega_2$.\\
For an infinite cardinal $\kappa$ let $H(\kappa)$ be the class of
sets hereditarily of cardinality $< \kappa$, i.e. $H(\kappa)=\{ X:
|TC(X)| < \kappa \}$ where $TC(X)$ is the transitive closure of
the set $X$. As in \cite{Sh2} there is an elementary submodel $N
\prec (H(\aleph_2),\varepsilon)$ such that
\begin{itemize}
\item $|N|=\aleph_0$;
\item $p, f_{\tilde{\bar{\nu}}, \tilde{\bar{c}}} \in N_0$;
\item $N=\bigcup_{n \in \omega} N_n$ with $N_n \prec
(H(\aleph_2),\varepsilon)$ elementary submodels such that $N_n \in
N_{n+1}$.
\end{itemize}
We let $\delta=N \cap \omega_1 \in S$ and $\delta_n = N_n \cap
\omega_1$ for $n \in \omega$. Choose $\eta_{\delta}$ such that $
\eta_{\delta}(n) \in [\delta_n, \delta_{n+1}]$
for all $n \in \omega$ and $\eta_{\delta}$ is simple and special.\\
As in \cite[Lemma 1.8]{Sh1} and \cite[Theorem 2.1]{Sh2} (see also
\cite{EkSh}) we define
inductively a sequence of finite sets of conditions in the following way:\\
In stage $n$ let $\eta_{\delta}(n)=\gamma$. We have a finite tree
$\left< p_t^n : t \in T_n \right> \in N_{n+1}$ of conditions and
let $k_{n+1}^{\delta}=k_n^{\delta} + |max(T_n)| +1$. Moreover, if
$t \in \max(T_n)$, then
\[ p_t^n \text{ forces a value to }
f_{\tilde{\bar{\nu}}, \tilde{\bar{c}}} \restriction
G_{\tilde{\bar{\nu}} \restriction_{\gamma + \omega}}.
\] We now choose $\nu_{\delta}(k_n^{\delta}+i)$ in $[
\gamma, \gamma +\omega ]$ for $i < t_n^{\delta}$ so that
$\nu_{\delta}$ becomes special and
$\rd(\eta_{\delta})=\rd(\nu_{\delta})$. We have that, if $t \in
\max(T_n)$, then
\[ p_t^n \Vdash
"f_{\tilde{\bar{\nu}}, \tilde{\bar{c}}}
(x_{\nu_{\delta}(k_n^{\delta}+i)}) -
\hat{x}_{\nu_{\delta}(k_n^{\delta}+i)}=\hat{w}\tilde{b}_{t,n,i}"
\]
for every $i < t_n^{\delta}$. By linear algebra we may choose
$a_l^{\delta,n}$ for $l < t_n^{\delta}$ such that
$\gcd(a_l^{\delta,n} : l < t_n^{\delta})=1$ and if $t \in
\max(T_n)$, then
\[ p_t^n \Vdash "\sum\limits_{i < t_n^{\delta}}a_i^{\delta,n} \tilde{b}_{t,n,i} =0". \]
Finally, we choose $c_{\delta}(n)$ arbitrary.  In the inverse
limit we hence obtain a triple $\left< \bar{\eta}', \bar{\nu}',
\bar{c}' \right>$ which we may increase to $\left<
\bar{\eta}'\left< \eta_{\delta} \right>, \bar{\nu}'\left<
\nu_{\delta} \right>, \bar{c}'\left< c_{\delta} \right> \right>$.
Now we can find $p^* \in P_{\alpha}$ above $p_t^n$ for some $t \in
\max(T_n)$ and all $n \in \omega$. Note that the $c_{\delta}$ was chosen arbitrary, so
there are $2^{\aleph_0}$ possible choices for the same $\nu_{\delta}$. \\
Now assume that $G$ is a generic filter containing the condition
$p^{*}$. Then $\hat{z}_{\delta,0} - f_{\tilde{\bar{\nu}},
\tilde{\bar{c}}}(z_{\delta,0}) \in \hat{w} \Z$. Moreover,
\[ n!f_{\bar{\nu},
\bar{c}}(z_{\delta,n+1})=f_{\bar{\nu}, \bar{c}}(z_{\delta,n}) +
\sum\limits_{i < t_n^{\delta}} a_i^{\delta,n} f_{\bar{\nu},
\bar{c}}(x_{\nu_{\delta}(k_n^{\delta} + i)}) \] for every $n \in
\omega$. Similarly, we have
\[ n!\hat{z}_{\delta,n+1}=\hat{z}_{\delta,n} + \sum\limits_{i < t_n^{\delta}}
a_i^{\delta,n} \hat{x}_{\nu_{\delta}(k_n^{\delta} + i)} +
c_{\delta}(n)\hat{w}.
\] Subtracting the two equations yields
\[ n!(f_{\bar{\nu},
\bar{c}}(z_{\delta,n+1})-\hat{z}_{\delta,n+1})=(f_{\bar{\nu},
\bar{c}}(z_{\delta,n}) - \hat{z}_{\delta,n}) + \sum\limits_{i <
t_n^{\delta}} a_i^{\delta,n} (f_{\bar{\nu},
\bar{c}}(x_{\nu_{\delta}(k_n^{\delta} +
i)})-\hat{x}_{\nu_{\delta}(k_n^{\delta} + i)}) -
c_{\delta}(n)\hat{w}.
\] But by our choice we have \[ \sum\limits_{i < t_n^{\delta}}
a_i^{\delta,n} (f_{\bar{\nu},
\bar{c}}(x_{\nu_{\delta}(k_n^{\delta} +
i)})-\hat{x}_{\nu_{\delta}(k_n^{\delta} + i)})= \sum\limits_{i <
t_n^{\delta}} a_i^{\delta,n} b_{t,n,i} =0. \] Therefore, we get
\begin{equation}
\label{equ2.1}
n!(f_{\bar{\nu},
\bar{c}}(z_{\delta,n+1})-\hat{z}_{\delta,n+1})=(f_{\bar{\nu},
\bar{c}}(z_{\delta,n}) - \hat{z}_{\delta,n}) -
c_{\delta}(n)\hat{w} \in \Z \hat{w}.
\end{equation}
Since $\Z$ is countable there exist generic filters $G_1$ and
$G_2$ (and corresponding triples $(\bar{\eta}^1$, $\bar{\nu}^1,
\bar{c}^1)$ and $(\bar{\eta}^2$, $\bar{\nu}^2, \bar{c}^2)$) such
that $c_{\delta}^1 \not= c_{\delta}^2$,
$\nu_{\delta}^1=\nu_{\delta}^2=\nu_{\delta}$ but
\[ \hat{z}_{\delta,0} - f_{\bar{\nu}^1,
\bar{c}^1}(z_{\delta,0}) = \hat{z}_{\delta,0} - f_{\bar{\nu}^2,
\bar{c}^2}(z_{\delta,0}) \in \Z\hat{w}. \] Let $n$ be minimal such
that $c_{\delta}^1(n) \not= c_{\delta}^2(n)$. Then an easy
induction using equation $(\ref{equ2.1})$ shows that

\[ f_{\bar{\nu}^1,
\bar{c}^1}(z^1_{\delta,l})- \hat{z}^1_{\delta,l}=f_{\bar{\nu}^2,
\bar{c}^2}(z^2_{\delta,l}) - \hat{z}^2_{\delta,l} \in \Z \hat{w}\]
for every $l \leq n$. Note that $\hat{z}_{\delta,i}$ depends on
$G_i$ $(i=1,2)$. We finally calculate
\[ n!(f_{\bar{\nu}^1,
\bar{c}^1}(z^1_{\delta,n+1})- \hat{z}^1_{\delta,n+1}) - n!(
f_{\bar{\nu}^2, \bar{c}^2}(z^2_{\delta,n+1}) -
\hat{z}^2_{\delta,n+1})\]
\[ = (f_{\bar{\nu}^1,
\bar{c}^1}(z^1_{\delta,n})- \hat{z}^1_{\delta,n}) -
(f_{\bar{\nu}^2, \bar{c}^2}(z^2_{\delta,n}) -
\hat{z}^2_{\delta,n}) + (c_{\delta}^1(n)-c_{\delta}^2(n))\hat{w}.
\] By equation $(\ref{equ2.1})$ we conclude \[ n!(f_{\bar{\nu}^1,
\bar{c}^1}(z^1_{\delta,n+1})- \hat{z}^1_{\delta,n+1}) - n!(
f_{\bar{\nu}^2, \bar{c}^2}(z^2_{\delta,n+1}) -
\hat{z}^2_{\delta,n+1})=
(c_{\delta}^1(n)-c_{\delta}^2(n))\hat{w}.\] However, the left side
is divisible by $2$ but the right side is $\hat{w}$ or $-\hat{w}$,
hence not divisible by $2$ - a contradiction. Note that all the differences are elements of the pure subgroup $\hat{w}\Z$ by equation $(\ref{equ2.1})$. Hence the above calculations take place in $\hat{w}\Z$ which is in the ground model, although the elements we are talking about come from different (incompatible) extension models. Thus the sequence
$(E)$ cannot be forced to split and this finishes the proof. \qed

\begin{corollary}
It is consistent with $ZFC$ and $GCH$ that there exist two
filtration-equivalent $\aleph_1$-separable abelian groups of
cardinality $\aleph_1$ such that one is Whitehead and the other
is not.
\end{corollary}

\proof Applying Theorem \ref{maintheorem} to a model $V$ of $GCH$
there is an extension model of $V$ in which there exist abelian
groups $G_{\bar{\eta}}$ and $G_{\bar{\nu}}$ for $\bar{\eta},
\bar{\nu} \in E$ such that $G_{\bar{\eta}}$ is a Whitehead group
but $G_{\bar{\nu}}$ is not. Since $\rd(\bar{\eta})=\rd(\bar{\nu})$
we deduce that $G_{\bar{\eta}}$ and $G_{\bar{\nu}}$ are filtration-equivalent by Proposition \ref{mainprop}. \qed

\goodbreak

\bigskip

\end{document}